\documentclass[12pt,notitlepage,a4paper,reqno,ps]{amsart}
\usepackage[utf8]{inputenc}
\usepackage[T1]{fontenc}
\usepackage[english]{babel}
\usepackage{amsmath}
\usepackage{amsfonts}
\usepackage{amssymb}
\usepackage{amsthm}
\usepackage{mathrsfs}
\usepackage{graphicx}
\usepackage{braket}
\usepackage{color}
\usepackage{amscd,epsfig,esint,graphics}
\usepackage{wrapfig}
\usepackage{verbatim}
\usepackage{mathtools}

\theoremstyle{plain}
\newtheorem{theorem}{Theorem}[section]
\newtheorem{lemma}[theorem]{Lemma}
\newtheorem{proposition}[theorem]{Proposition}
\newtheorem{example}[theorem]{Example}

\newtheorem{remark}[theorem]{Remark}
\newtheorem{definition}[theorem]{Definition}

\theoremstyle{definition}
\theoremstyle{remark}
\numberwithin{equation}{section}

\newcommand{\R}{\mathbb{R}}
\newcommand{\N}{\mathbb{N}}
\newcommand{\K}{\mathbb{K}_\psi(\Omega)}

\newcommand{\e}{\varepsilon}
\newcommand{\C}{\mathcal{C}}
\newcommand{\om}{\Omega}

\newcommand{\dive}{\textnormal{div}}
\newcommand{\supp}{\textnormal{supp}\,}

\newenvironment{samuelerev}{\color{cyan}}{\color{black}}
\newcommand{\bsamr}{\begin{samuelerev}}
\newcommand{\esamr}{\end{samuelerev}}

\newenvironment{todorev}{\color{magenta}}{\color{black}}
\newcommand{\btodo}{\begin{todorev}}
\newcommand{\etodo}{\end{todorev}}

\definecolor{dg}{rgb}{0.01, 0.75, 0.24}

\title[Necessary condition for autonomous obstacle problems with general growth]{A necessary condition for extremality \\ of solutions to autonomous obstacle \\ problems with general growth}

\author[S.\,Ricc\`o]{Samuele Ricc\`o}
\address[S.\,Ricc\`o]{TU Wien, Institute of Analysis and Scientific Computing, Wiedner Hauptstraße 8-10, 1040, Vienna, Austria.}
\email{samuele.ricco@tuwien.ac.at}

\author[A.\,Torricelli]{Andrea Torricelli}
\address[A.\,Torricelli]{Universit\`a degli Studi di Modena e Reggio Emilia, Dipartimento di Scienze Fisiche, Informatiche e Matematiche, via Campi 213/b, 41125, Modena, Italy.}
\email{andrea.torricelli@unimore.it}

\begin{document}
	
\baselineskip 3.4ex
\vspace{0.5cm}

\maketitle
	
\begin{abstract}
Let us consider the autonomous obstacle problem
\begin{equation*}
\min_v \int_\Omega F(Dv(x)) \, dx
\end{equation*}
on a specific class of admissible functions, where we suppose the Lagrangian satisfies proper hypotheses of convexity and superlinearity at infinity. Our aim is to find a necessary condition for the extremality of the solution, which exists and it is unique, thanks to a primal-dual formulation of the problem. The proof is based on classical arguments of Convex Analysis and on Calculus of Variations' techniques.
\end{abstract}


\section{Introduction}
\label{intro}

\noindent
In this manuscript we consider autonomous variational obstacle problems of the form
\begin{equation*}
\min \left\{ \int_\Omega F(Dv(x)) \, dx \, : \, v \in \K \right\},
\end{equation*}
where $\Omega$ is an open and bounded subset of $\R^n$ and $\K$ is the class of admissible functions, defined as
\begin{equation}
\label{class} \tag{K}
\K = \left\{ v \in W_{u_0}^{1,1}(\Omega) \, : \, v \ge \psi \, \text{ a.e. on } \Omega, \, F(Dv) \in L^1(\Omega) \right\},
\end{equation}
where $u_0 \in W^{1,1}(\Omega)$ is a boundary datum such that $F(Du_0) \in L^1(\Omega)$ and where $\psi \in W^{1,1}(\Omega)$ is a function called {\it obstacle} such that $F(D\psi) \in L^1(\Omega)$. The main focus we have in this paper is how minimizers of the above-defined constrained problem could be characterized, exploiting a primal-dual formulation of the optimization problem. We suppose the Lagrangian satisfies a superlinear growth condition at infinity, although it is not subjected to any growth condition from above. Moreover, we assume that the Lagrangian is bounded from below by a given convex and superlinear function, that brings the Lagrangian to inherit the same degree of convexity that it has.
\\
There are many works about regularity theory in variational problems and elliptic systems with non-standard growth, but the papers which paved the way were the famous \cite{M89} and \cite{M91} by Marcellini; since they were published, a lot of new ideas have been applied to this research branch and many results have been proved in several directions (see for example \cite{M96} and \cite{M2021} by Marcellini, or \cite{Ming} by Mingione for a general exposition and further references). However, regarding the obstacle problems there are still some issues which have not been studied in an exhaustive way yet. One of these issues deals with the relation between minima and extremals: it is common knowledge that, for both the constrained and unconstrained problems, the regularity of the solutions often comes from the fact that the minimizers are extremals too, i.e. they solve a corresponding variational equality or inequality. Note though that there are examples of variational problems whose minimizers do not satisfy the Euler-Lagrange equation in the weak sense, as proved by Ball and Mizel in \cite{BM}. While in the case of standard growth conditions the situation is well established (see for instance the book \cite{Dacorogna} by Dacorogna), in the case of non-standard growth conditions, the relation between extremals and minima is an issue that requires a careful investigation.
\\
In 2014 Carozza, Kristensen and Passarelli di Napoli investigated exactly this topic in \cite{CKPdN2014} in the case of convex integral functionals, with the aim of showing that their minimizers are characterized to be the energy solutions to the Euler-Lagrange systems for the functionals under non-standard growth conditions. The main tool they use is a particular regularization procedure: the integrand $F$ is approximated by a sequence of strictly convex and uniformly elliptic integrands $F_k$ which satisfy standard $p-$growth conditions and whose minimizers $u_k$ strongly converge to the minimizer $u$ in $W^{1,p}$. With that said, according to the standard duality theory for convex problems, every such minimizer $u_k$ is associated to a row-wise solenoidal matrix field denoted by $\sigma_k$. Finally, for the pairing $(Du_k, \sigma_k)$, suitable pointwise estimates that are preserved while passing to the limit are proved. Such estimates then provide conditions which allow the Euler-Lagrange system to hold for an $F$-minimizer. In a subsequent paper, i.e. \cite{CKPdN2015}, the same achievement has been carried on under more general growth assumptions, covering a wide class of functionals, from those with almost-linear growth to the ones with exponential growth and beyond, by the use of Ekeland variational principle and Young measures, to obtain the necessary estimates for the pairing $(Du_k, \sigma_k)$ to be able to pass to the limit.
\\
In \cite{CKPdN2014}, \cite{CKPdN2015}, by Carrozza, Kristensen and Passarelli di Napoli, and \cite{EPdN}, by Eleuteri and Passarelli di Napoli, the concept of convex duality is exploited. While this seems to be very natural, its use is not so common in the context of convex variational integrals with non-standard growth conditions. We should note though that, before these papers, also in \cite{BC} and \cite{BCM} the authors (respectively Bonfanti and Cellina for the first one and Bonfanti, Cellina and Mazzola for the second) make use of it and, in particular, in \cite{BCM} is also addressed the question of energy-extremality of minimizers. They work in the context of more general variational integrals $v \, \longmapsto \, \int_\Omega F(x, v, Dv) \, dx$ in the multi-dimensional scalar case with $n \ge 2$, $N = 1$ under convexity and regularity hypotheses on the gradient. Talking about more general functionals, we remark that the results obtained in \cite{CKPdN2014} can be generalized to minimizers of the general autonomous convex variational integral $\int_\Omega F(v, Dv) \, dx$, just under the hypothesis that the integrand $F = F(\eta, \, \xi)$ is jointly convex. Similar remarks, together with the precise statements and sketches of the proofs, were given in \cite{CKPdN11} by the same authors.
\\
Eventually, in the paper \cite{EPdN} Eleuteri and Passarelli di Napoli address the analogue issue of \cite{CKPdN2014} in the case of constrained minimizers with a very general obstacle quasi-continuous up to a subset of zero capacity. Let us mention that relying on techniques of convex analysis, Scheven and Schmidt in \cite{SS16} and \cite{SS} investigated the Dirichlet minimization problem for the total variation and the area functional with one-sided obstacle. The main point is that they were able to identify certain dual maximization problems for bounded divergence-measure fields and to establish duality formulas and point-wise relations between (generalized) $BV$ minimizers and dual maximizers. Their results are very general and apply to very general obstacles, such as $BV$ obstacles and the obstacle considered in \cite{EPdN}; the proofs of their results crucially depend on a new version of the Anzellotti-type pairing (see \cite{A84}) which involves general divergence measure fields and specific representatives of $BV$ functions, by employing several fine results on capacities and one-sided approximation. This framework is proved to be the right one in order to extend the results in \cite{CKPdN2014} to very general obstacle problems, as long as, by means of the Anzellotti-type pairing, they are able to express the natural counterpart of the variational inequality in this very general setting, which will reduce to the usual one once they have the right summability for the functions involved.
\\
We were inspired by them in order to try to extend the results of \cite{CKPdN2015} in the constrained optimization problem, but we make use of different hypotheses on the obstacle and on the Lagrangian, in particular a superlinear growth condition at infinity and convexity guaranteed by the hypotheses on the function that bounds the Lagrangian from below. As in \cite{CKPdN2014}, \cite{CKPdN2015} and \cite{EPdN}, we use the concept of convex duality in various steps of the proof. The most challenging knots in the proof under our hypotheses are the passage to the limit, where we have to pay attention to the presence of the obstacle, and the proof of the variational inequality.
\\
Our paper is organized as follows: in Section \ref{main_prob} we state the hypotheses and the definitions we need in the following to prove our results and we state the two new theorems of the paper, together with some examples of functionals to which our first result applies (Subsection \ref{exam}). Section \ref{prelim} holds some preliminary results about Convex Analysis and Young measures, together with some notations. The proof of the main result is then contained in Section \ref{proof_main}, which is divided in six small parts in order to give a more clear stream of reading, while the proof of our second theorem is contained in Section \ref{dimostr-2}.
\noindent


\section{Statement of the Main Result}
\label{main_prob}

\noindent
Let us consider the problem
\begin{equation}
\label{oprob}
\min_{v \, \in \, \K} \, \int_\Omega F( Dv(x) ) \, dx,
\end{equation}
where $\Omega \subset \R^n$ is an open and bounded set, $F : \R^n \longrightarrow [0,+\infty)$ is a $\C^1$ function and $u_0 \in W^{1,1}(\Omega)$ is a boundary datum such that $F(Du_0) \in L^1(\Omega)$. Moreover, the function $\psi \in W^{1,1}(\Omega)$ is called {\it obstacle} and it is such that $F(D\psi) \in L^1(\Omega)$. The class of the admissible functions $\K$ is defined as in the introduction, but we remark it here:
\begin{equation}
\label{class} \tag{K}
\K \, := \, \left\{ v \in W_{u_0}^{1,1}(\Omega) \, : \, v \ge \psi \, \text{ a.e. on } \, \Omega, \, F(Dv) \in L^1(\Omega) \right\}.
\end{equation}
\noindent
We suppose that there exists a $\C^1$ and strictly convex function $\phi : \R^n \longrightarrow [0,+\infty)$ such that
\begin{equation}
\label{def-phi}
\phi(\xi) \, := \, \theta(|\xi|) \qquad \textnormal{for all} \quad \xi \in \R^n
\end{equation}
for a function $\theta : [0,+\infty) \longrightarrow [0,+\infty)$ superlinear at infinity, and we suppose that
\begin{equation}
	\label{H1} \tag{H1}
	F - \phi \qquad \text{is a convex function,}
\end{equation}
which implies that $F$ is strictly convex. Moreover, we suppose that there exists $c \in \R$ such that 
\begin{equation}
	\label{H1'} \tag{H2}
	F(\xi) \, \ge \,  \frac{1}{2}\phi(\xi) +c,
\end{equation}
for all $\xi \in \R^n$ such that $|\xi|$ is large enough. Clearly, this hypothesis implies that $F$ is superlinear at infinity.
\begin{remark}
We point out that hypothesis \eqref{H1'} can be immediately replaced by the superlinearity at infintiy of $F$. We decided to use this formulation though in order to highlight the fact that $F$ inherits the same properties of $\phi$.
\end{remark}
\noindent
Now we define the space
\begin{equation*}
S_{-}(\Omega) \, := \, \{ \sigma: \Omega \longrightarrow \R^n \, : \, \dive \, \sigma \le 0 \ \text{ in distributional sense} \}.
\end{equation*}
Fixed $U \in W^{1,1}(\Omega)$ and $\sigma \in S_{-}(\Omega)$, we also define the measure $[\![\sigma, U]\!]_{u_0} (\overline{\Omega})$ on $\overline{\Omega}$ as
\begin{equation}
\label{misura-sigma}
[\![\sigma, U]\!]_{u_0} (\overline{\Omega}) \, := \, \int_\Omega (U - u_0) \, d (- \dive \, \sigma) + \int_\Omega \braket{\sigma, \, Du_0} \, dx.
\end{equation}
\begin{remark}
It is worth noticing that, in general, $-\dive \sigma$ is a measure and that the measure $[\![\sigma, U]\!]_{u_0} (\overline{\Omega})$ may take the value $+\infty$ for certain choices of $U \in W^{1,1}(\Omega)$ and $\sigma \in S_{-}(\Omega)$. Moroever, the measure \eqref{misura-sigma} can actually be defined for a general $u_0\in W^{1,1}(\Omega)$ with no further hypotheses.
\end{remark}
\noindent
Chosen $U \in W^{1,1}_{u_0}(\Omega)$, then $[\![\sigma, U]\!]_{u_0} (\overline{\Omega})$ is equal to $\braket{\sigma, \, DU} \in L^1(\Omega)$ since, by definition of distributional divergence,
\begin{equation}
\label{parts}
\int_\Omega \varphi \, d(- \dive \, \sigma) \, = \, \int_\Omega \braket{\sigma, \, D\varphi} \, dx \qquad \forall \, \varphi \in W_0^{1,1}(\Omega).
\end{equation}
Finally we observe that under the assumptions \eqref{H1} and \eqref{H1'}, the functional
\begin{equation*}
I(v) \, := \, \int_\Omega F( Dv(x) ) \, dx
\end{equation*}
 is proper, convex and lower semicontinuous on $W^{1,1}(\Omega)$ and thus, given $u_0 \in W^{1,1}(\Omega)$ such that $F(Du_0)\in L^1(\Omega)$, the existence and uniqueness of the minimizer $u$ in the convex space $\K$ are granted. A necessary condition for the extremality of this minimizer is stated in the following Theorem, which is our main result.

\begin{theorem}
\label{mtheorem}
Let $F$ be a non-negative function in $\C^1 (\R^n)$, satisfying \eqref{H1} and \eqref{H1'} with $\phi$ defined as in \eqref{def-phi}, and let $u_0 \in W^{1,1}(\Omega)$ be such that $F(Du_0), F(t \, Du_0) \in L^1(\Omega)$ for some $t > 1$. Then, for the unique minimizer $u \in W_{u_0}^{1,1}(\Omega)$ of the minimum problem \eqref{oprob} it holds
\begin{equation}
\label{dimdim}
F^*( F' (Du) ) \in L^1(\Omega), \qquad \braket{F'(Du), Du} \in L^1(\Omega)
\end{equation}
and
\begin{equation}
\label{div_neg}
\dive \, F'(Du) \, \le \, 0
\end{equation}
in distributional sense. Moreover it holds the following identity
\begin{equation}
\label{min_max_finale}
\int_\Omega F(Du) \, dx \, = \, [\![F'(Du), \psi]\!]_{u_0} (\overline{\Omega}) - \int_\Omega F^*(F'(Du)) \, dx.
\end{equation}
\end{theorem}

\noindent
Since superlinearity and Lipschitz continuity cannot coexist as properties of the same function, it is only natural to wonder whether the same results hold true if we consider the Lipschitz continuity of the integrand instead of the superlinearity at infinity. To answer this question we prove the following theorem. It is important to notice that the Lipschitz continuity (and thus lack of superlinearity at infinity) does not ensure the existence of the minimizer and so we have to assume its existence (and uniqueness).

\begin{theorem}
\label{dualproblem}
Let $G \in \C^1(\R^n)$ be strictly convex and Lipschitz continous and such that there exists a solution $u$ of the problem
\begin{equation}
\label{gprob}
\min_{v \, \in \, \K} \, \int_\Omega G( Dv(x) ) \, dx,
\end{equation}
where the class $\K$ is defined as in \eqref{class}. Let $u_0\in W^{1,1}(\Omega)$ be such that $G(Du_0)\in L^1(\Omega)$. Then the following holds true
\begin{equation}
	\label{prob-dual}
	\min_{v \in \K} \int_{\Omega} G(Dv) \, dx = \max_{\sigma \in S_{-}(\Omega)} \left ([\![\sigma, {\psi}]\!]_{u_0}(\overline{\Omega}) - \int_{\Omega} G^*(\sigma) \, dx \right ).
\end{equation}
In particular,
\begin{equation}
	\label{prob-dual2}
	\int_{\Omega} G(Du) \, dx = [\![G'(Du), {\psi}]\!]_{u_0}(\overline{\Omega}) - \int_{\Omega} G^*(G'(Du)) \, dx. 
\end{equation}
\end{theorem}
\noindent
We remark that if $G$ is Lipschitz continuous there is no need to prove the conditions \eqref{dimdim} and \eqref{div_neg}. Indeed, the first follows from the fact that if $G$ is Lipschitz continuous then $G'\in L^\infty(\Omega)$ and thus
\begin{align*}
	\braket{G'(Du), Du} \in L^1(\Omega),
\end{align*}
considered together with \eqref{fineq2}, that holds since $G$ is convex and lower semicontinuous. The second, namely \eqref{div_neg}, instead follows from \eqref{varineq}.


\subsection{Some examples of functionals}
\label{exam}

As stated in the introduction, namely Section \ref{intro}, the general growth assumptions we are considering on the convex Lagrangian $F$, in particular \eqref{H1'}, allow us to cover a wide class of functionals from those with almost-linear growth to the ones with exponential growth and beyond. For specific examples of functionals with general growth, see Section 3 of \cite{EMMP}. We will now state some examples of Lagrangians to which Theorem \ref{mtheorem} applies to. Note that these functions will be always called $F_i$ for a certain index $i \in \N$, while auxiliary functions will be denoted with different letters.

\begin{example}
For every $\xi \in \R^n$, we can consider
\begin{equation*}
F_1(\xi) := |\xi|^\alpha \qquad \textnormal{ with } \ \alpha > 1.
\end{equation*}
\end{example}

\begin{example}
For every $\xi \in \R^n$, we can consider
\begin{equation*}
F_2(\xi) := \cosh|\xi|.
\end{equation*}
\end{example}

\begin{example}
For every $\xi \in \R^n$ we consider the function
\begin{equation*}
G(\xi) := |\xi| \ln |\xi|
\end{equation*}
and we refer to one of its global minimum points as $\xi_0$, namely such that $|\xi_0| = \frac{1}{e}$ and $G(\xi_0) = -\frac{1}{e}$. If we define 
\begin{equation*}
F_3(\xi) := (|\xi| + |\xi_0|) \ln (|\xi| + |\xi_0|) - G(\xi_0) = \left(|\xi| + \frac{1}{e}\right) \ln \left(|\xi| + \frac{1}{e}\right) + \frac{1}{e},
\end{equation*}
then it satisfies the hypotheses of Theorem \ref{mtheorem}.
\end{example}


\section{Preliminaries}
\label{prelim}
\noindent
We denote with $\C(\R^n)$ the space of the continuous functions from $\R^n$ to $\R$ and with $\C_0(\R^n)$ the closure of $\C(\R^n)$ with respect to the maximum norm. We will also denote with $M(\R^n)$ the space of finite Radon measures on $\R^n$ and, given a generic measure $\mu$, we will denote its support with $\supp \mu$. We denote with $Pr(\R^n)$ the space of the probability measures defined on the Borel sets of $\R^n$. Moreover, we denote with $L^\infty_w(\Omega, M(\R^n))$ the set of the essentially bounded and weakly$^*$ measurable functions from $\Omega$ to $M(\R^n)$. Finally, we denote with $\overline{\R}$ the extended real number line.

\subsection{Convex analysis}
In this subsection we will state some known results about convex functions and conjugate functions. For more details see \cite{ET}.
%
%
\begin{definition}
Given a continuous function $F : \R^n \longrightarrow \, \overline{\R}$, its polar function is the function $F^* : \R^n \longrightarrow \, \overline{\R}$ defined as
	\begin{equation*}
		F^*(z) := \sup_{\xi \, \in \, \R^n} \left[ \braket{z, \xi} - F(\xi) \right] \qquad \forall \, z \in \R^n.
	\end{equation*}
The polar function $F^{**} : \R^n \longrightarrow \, \overline{\R}$ of $F^*$ is called bipolar function of $F$ and it is defined as
\begin{equation*}
	F^{**}(\xi) := \sup_{z \, \in \, \R^n} \left[ \braket{\xi, z} - F^*(z) \right] \qquad \forall \, \xi \in \R^n.
\end{equation*}
\end{definition}

\begin{remark}
It is well known that $F^*$ is always a convex and lower semicontinuous function. Moreover, it is possible to prove that $F^{**} = F$ if and only if $F$ is convex and lower semicontinuous for every $\xi \in \R^n$.
\end{remark}

\noindent
By the definitions of polar and bipolar functions, it is clear that for every $\eta, \, \xi \in \R^n$, it also holds
\begin{equation}
\label{fineq}
\braket{\xi, \eta} \le F^*(\eta) + F^{**}(\xi).
\end{equation}
This inequality is known as Fenchel inequality. The equality holds true when $\eta \in \partial F^{**}(\xi)$, where $\partial F^{**}(\xi)$ is the subgradient of $F^{**}$ in $\xi$ (see again \cite{ET} for details). In particular, if $F\in \C^1(\R^n)$ and it is convex, then in \eqref{fineq} the equality holds true for $\eta = F'(\xi)$, i.e.
\begin{equation}
\label{fineq2}
\braket{\xi, F'(\xi)} = F^*(F'(\xi)) + F^{**}(\xi).
\end{equation}
\noindent
The polar function also have an additional important property that we will need in the following, in particular given by the following result, which is a corollary of \cite[Lemma 3.1]{GT}.
\begin{lemma}
\label{superlinear}
If $F : \R^n \longrightarrow \, \overline{\R}$ is convex, non-negative, superlinear and $F(0)=0$ and has as effective domain the whole $\R^n$, then its polar $F^*$ verifies the same properties of $F$.
\end{lemma}
\begin{remark}
It is worth noticing that the hypothesis $F(0)=0$ is not necessary for the proof of the Lemma, in the sense that it is only a property the polar function inherits from the original function $F$.
\end{remark}
\noindent
We conclude this subsection stating the Ekeland variational principle. The original Theorem is proved for Banach spaces and can be found in \cite[Corollary 6.1]{ET}. Our space $\K$ is not a Banach space (it is not even a linear space), but if endowed with the distance induced by the norm
\begin{align*}
	\|v\|_{W^{1,1}(\Omega)} := \int_\Omega |Dv| \, dx,
\end{align*}
it is instead a complete metric space, indeed it is a closed subset of a complete metric space. This follows from Poincaré inequality and the fact that $\K\subset W^{1,1}_{u_0}(\Omega)$. Moreover, since $\K$ is closed and convex it is also weakly closed. The result that follows is a version of Ekeland variational principle stated for metric spaces, whose proof can be found in \cite[Section 1.4]{Z}.
\begin{theorem}[Ekeland Variational Principle]
\label{Ekeland}
Let $(V, d)$ be a complete metric space and $F : V \longrightarrow \R$ be a lower semicontinuous function bounded from below. Given $\e > 0$ and $v \in V$ such that 
	\begin{align*}
		F(v) \, \le \, \inf_V F + \e,
	\end{align*}
then for every $\lambda > 0$ there exists $v_\lambda \in V$ such that
\begin{itemize}
	\item[i)] $d(v,v_\lambda) \le \lambda$,
	\item[ii)] $F(v_\lambda) \le F(v)$,
	\item[iii)] $v_\lambda$ is the unique minimizer of the functional $v \longmapsto F(v) + \e \lambda^{-1} d(v, v_\lambda)$.
\end{itemize}
\end{theorem}
\noindent

\subsection{Generalized Young measures}
In order to keep this paper self-contained, we now state some known results about Generalized Young measures. For more details on Young measures and Generalized Young measures see \cite{R}. Note that in this section (and everytime we work with generalized Young measures) the braket $\braket{\cdot,\cdot}$ will denote the duality product instead of the usual scalar product in $\R^n$.

\begin{definition}
	Given a Carathéodory function $f:\R^n\to\R$ with linear growth we call Recession function (or Recession integrand) the function $f^\infty$ defined as 
	\begin{align*}
		f^\infty(\xi):=\lim_{t \to +\infty}\frac{f(t \, \xi)}{t},
	\end{align*}
if this limit exists.
\end{definition}

\begin{definition}
	A generalized Young measure on the open and bounded set $\Omega\subset\R^n$ is a triple $\nu=((\nu_x)_{x\in\om}, \lambda, (\nu_x^\infty)_{x\in\om})$ where $\nu_x\in Pr(\R^n)$ for a.e. $x\in \Omega$, $\lambda$ is a positive finite measure on $\overline{\Omega}$ and $\nu_x^\infty\in Pr(\mathbb{S}^{n-1})$ such that
	\begin{itemize}
		\item[i)] $x\mapsto \nu_x$ is weakly$^*$ measurable with respect to Lebesgue measure,
		\item[ii)] $x\mapsto \nu_x^\infty$ is weakly$^*$ measurable with respect to $\lambda$,
		\item[iii)] $x\mapsto \braket{|\cdot|,\nu_x}\in L^1(\Omega)$.
	\end{itemize}
\end{definition}
\noindent
The next Proposition includes some well known results about generalized Young measures.

\begin{proposition}
	\label{minestrone}
	Given a sequence $(v_k)_k\subset W^{1,1}(\Omega)$ such that $\sup_{k}||Dv_k||_{L^1(\Omega)}<+\infty$, then there exist a non relabeled subsequence $(v_k)_k$ and a generalized Young measure $\nu=((\nu_x)_{x\in\om}, \lambda, (\nu_x^\infty)_{x\in\om})$ such that it holds
	\begin{align*}
		\lim_{k \to +\infty}\int_\Omega f(Dv_k)dx=\int_\Omega\braket{f,\nu_x}dx+\int_{\overline{\Omega}} \braket{f^\infty,\nu_x^\infty}d\lambda,
	\end{align*}
for any function with linear growth, i.e. such that there exists $m>0$ for which $|f(\xi)|\le m(1+|\xi|)$. Moreover $(v_k)_k$ is equi-integrable if and only if $\lambda=0.$
\end{proposition}

\begin{remark}
	\label{Ym}
	Given $(v_k)_k\subset W^{1,1}(\Omega)$ such that $\sup_{k}\|Dv_k\|_{L^1(\Omega)}<+\infty$, then the sequence $(v_k)_k$ is bounded in $BV$ and we can extract a non relabeled subsequence $(v_k)_k$ converging weakly$^*$ in $BV$ to a function $v$. Moreover, if $\lambda=0$ then the bycenter of the generalized Young measure $\nu$ generated by the previous proposition is 
	\begin{align*}
		[\nu](x):=[\nu_x]=\braket{Id,\nu_x}=Dv
	\end{align*}
\end{remark}



\section{Proof of Theorem \ref{mtheorem}}
\label{proof_main}

\noindent
In this section we will prove our main result, namely Theorem \ref{mtheorem}. The proof is conveniently divided into six steps in order to let the stream of reading clearer. In the first steep we build a sequence of approximating functions, while in the second step we work with sequences of perturbed problems. In the third step we pass to the limit and in the last three steps we prove the three theses of the Main Theorem, i.e. respectively \eqref{dimdim}, \eqref{div_neg} and \eqref{min_max_finale}.

\subsection{Construction of approximating functions}
\label{quat_zero}
We consider the polar function $F^*$ of $F$. Thanks to Lemma \ref{superlinear}, since $F$ is superlinear, convex and $\C^1$, then its polar
\begin{equation*}
F^*(z) := \sup_{\xi \, \in \, \R^n} \left( \braket{z, \xi} - F(\xi) \right) \qquad \forall \, z \in \R^n
\end{equation*}
is a real-valued, strictly convex and superlinear at infinity function. Fixed $k \in \N$ and $\xi \in \R^n$, we can define
\begin{equation*}
\overline F_k^{**}(\xi) := \sup_{|z| \, \le \, k} \, \left( \braket{\xi, z} - F^*(z) \right).
\end{equation*}
We can observe that $\overline F^{**}_k$ is a real-valued, convex and $k$-Lipschitz function and, since $F^{**}$ is lower semicontinuous, then 
\begin{equation*}
\overline F_k^{**} \stackrel{k \to \infty}{\longrightarrow} F^{**} = F \quad \text{pointwise,}
\end{equation*}
sinec $F$ is lower semicontinuous and convex. Now we define 
\begin{equation}
	\label{Fkpositive}
	\overline{G}^{**}_k (\xi) := \max \left\{ \overline F^{**}_k(\xi), \theta(|\xi|) \right\}.
\end{equation}
Again, $\overline{G}^{**}_k \stackrel{k \to \infty}{\longrightarrow} F^{**} = F$ pointwise, but for each $k \in \N$ there also must exists $r_k > 0$ such that $r_k \stackrel{k \to \infty}{\longrightarrow} +\infty $ and such that
\begin{equation*}
	\overline{G}^{**}_k(\xi) = \theta(|\xi|)
\end{equation*}
when $|\xi| \ge r_k$. Now we define
\begin{equation*}
	H_k (\xi) :=
	\begin{dcases}
		\overline{G}^{**}_k(\xi) \qquad \qquad \,\, &\text{if  } |\xi| < r_k, \\
		\frac{\theta(r_k)}{r_k}\, |\xi|  \qquad &\text{if  } |\xi| \ge r_k.
	\end{dcases}
\end{equation*}
Again, it is possible to prove that $H_k$ is a convex and $m_k - $Lipschitz function, with
\begin{equation*}
	m_k := \frac{\theta(r_k)}{r_k}
\end{equation*}
for all $k \in \N$. Now we regularize $H_k$ by means of the convolution kernels
\begin{equation*}
\Phi_\e (\xi) := \e^{-n} \, \Phi \left( \frac{\xi}{\e} \right),
\end{equation*}
where
\begin{equation*}
\Phi (\xi) :=
\begin{dcases}
& c \, \exp \left( {\frac{1}{|\xi|^2 - 1}} \right) \qquad \text{if  } |\xi| < 1, \\
& 0 \qquad \qquad \qquad \qquad \quad \, \text{if  } |\xi| \ge 1 \\
\end{dcases}
\end{equation*}
and where $c$ is chosen such that $\int_{\R^n} \Phi(\xi) \, d\xi = 1$. In particular, we consider the function $\Phi_\e * H_k (\xi)$ and we remark that this is a convex, $\C^\infty$ and $m_k - $Lipschitz function for which it holds
\begin{equation}
\label{stime-acca}
H_k (\xi) \le \Phi_\e * H_k(\xi) \le H_k(\xi) + \e \, m_k,
\end{equation}
for each $k \in \N$ and $\xi \in \R^n$. If we define
\begin{align*}
& \delta_k := \frac{1}{k^2m_k}, \\
& \mu_k := \frac{1}{k - 1}
\end{align*}
and 
\begin{equation}
\label{effe-kappa}
F_k (\xi) := \Phi_{\delta_k} * H_k (\xi) - \mu_k,
\end{equation}
then it is possible to prove that it holds true that $F_k(\xi) \le F_{k+1}(\xi)$ for all $k \in \N$ and $\xi \in \R^n$ (see \cite{CKPdN2014}) and that $F_k\nearrow F$ pointwise as $k\to \infty$.


\subsection{Perturbed problems}
\label{quat_uno}

\noindent
Since, by construction, the functions $F_k$ are not superlinear, we can not guarantee the existence of solutions to the obstacle problems
\begin{equation*}
\min_{w \, \in \, \K} \int_\Omega F_k(Dw) \, dx.
\end{equation*}
In order to bypass this issue, we will use a version of Ekeland variational principal for metric spaces, namely Theorem \ref{Ekeland}. We define for every $k \in \N$ the integral functional
\begin{align*}
	I_k(w) := \int_\Omega F_k(Dw) \, dx, \qquad \textnormal{for all } w \in \K.
\end{align*}
For every $k\in\N$ we have that $F_k\le F$, so
\begin{align}
	\label{bound}
	\inf_{w\in\K} I_k(w)\le\inf_{w\in\K} I(w)=\int_\Omega F(Du) \, dx.
\end{align}
By definition of infimum, we can find for every $k \in \N$ a function $v_k \in \K$ such that
\begin{align*}
	I_k(v_k) \le \inf_{w\in\K} I_k(w)+\frac{1}{k}.
\end{align*}
Noticing that for every $k\in\N$ holds
\begin{align*}
	F_k(\xi)\ge\min\left\{ \frac{\theta(r_k)}{r_k}\, |\xi|, \theta(|\xi|) \right\} \quad \textnormal{for all } \ \xi \in \R^n,
\end{align*}
then by \eqref{bound} we get that $(v_k)_k\subset\K$ must be bounded in $W^{1,1}(\Omega)$, which implies that there exists a not relabeled subsequence $(v_k)_k$ converging weakly$^*$ in $BV$ to some $v$. Moreover, by Proposition \ref{minestrone} we get that there exists a generalized Young measure $\nu=((\nu_x)_{x\in\om}, \lambda, (\nu_x^\infty)_{x\in\om})$ such that, for every $k \in \N$,
\begin{align}
	\label{gym}
	  \lim_{j \to +\infty}\int_\Omega F_j(Dv_j)\,dx \ge\lim_{j \to +\infty}\int_\Omega F_k(Dv_j)\,dx=\int_\Omega\braket{F_k,\nu_x}\,dx+\int_{\overline{\Omega}} \braket{F_k^\infty,\nu_x^\infty}\,d\lambda.
\end{align}
We observe now that, by construction,
\begin{align*}
F_k^\infty(\xi)=\lim_{t \to +\infty}\frac{F_k(t \, \xi)}{t} & = \frac{\theta(r_k)}{r_k}|\xi|,\\
\braket{F_k^\infty,\nu_x^\infty} & = \frac{\theta(r_k)}{r_k}.
\end{align*}
Moreover,
\begin{align}
	\label{acaso1}
	\lim_{j \to +\infty}\int_\Omega F_j(Dv_j)\,dx\le \lim_{j \to +\infty} \left(\inf_{w\in\K} I_j(w)+\frac{1}{j} \right) \le \int_\Omega F(Du)\,dx,
\end{align}
so passing to the limit as $k\to\infty$ in \eqref{gym}, by the Monotone Convergence Theorem and the fact that $\frac{\theta(r_k)}{r_k} \to +\infty$ as $k \to +\infty$, we get that $\lambda=0$ and 
\begin{align}
	\label{acaso2}
	\int_\Omega\braket{F,\nu_x}dx\le \lim_{j \to +\infty} \left(\inf_{w\in\K} \int_\Omega F_j(Dw)\,dx\right).
\end{align}
Since $\lambda=0$, then $(Dv_k)_k$ is equi-integrable by Proposition \ref{minestrone}, implying that $(v_k)_k$ is converging weakly in $W^{1,1}(\Omega)$ to $v$. Since $\K$ is weakly closed then also $v\in\K$. Furthermore, by Jensen's inequality and Remark \ref{Ym} we get
 \begin{align}
 	\label{lbound}
 	\int_\Omega F(Dv) \, dx \le \int_\Omega\braket{F,\nu_x} \, dx.
 \end{align}
This proves that
\begin{itemize}
	\item[i)] $v=u$,
	\item[ii)] $\inf_{w\in\K}\int_\Omega F_j(Dw) \, dx \longrightarrow \int_\Omega F(Du) \, dx$ as $j \to +\infty$.
\end{itemize}
Indeed by \eqref{acaso1}, \eqref{acaso2}, and \eqref{lbound} we have that 
\begin{align*}
\int_\Omega F(Dv)\,dx\le\int_\Omega\braket{F,\nu_x}dx\le \lim_{j \to +\infty} \left(\inf_{w\in\K} \int_\Omega F_j(Dw)\,dx\right)\le\int_\Omega F(Du)\,dx,
\end{align*}
and since $u$ is a minimizer and it is unique we get $(i)$ and $(ii)$. On the other hand, $(ii)$ allows us to say that there exists $(\e_k)_k\subset \R$ such that $\e_k\to 0$ as $k \to \infty$ and
\begin{align*}
	\inf_{w\in\K}\int_\Omega F_k(Dw)\,dx\le \int_\Omega F_k(Du)\,dx\le \int_\Omega F(Du)\,dx=\e_k^2+\inf_{w\in\K}\int_\Omega F_k(Dw)\,dx.
\end{align*}
By applying Theorem \ref{Ekeland} with $\lambda=\e_k$, we get that there exists a sequence $(u_k)_k \subset \K$ such that $D u_k\to D u$ in $L^1(\Omega)$ and such that
\begin{align*}
	\int_\Omega F_k(Du_k)\,dx\le\int_\Omega F_k(Du)\,dx \le \int_\Omega F(Du)\,dx.
\end{align*}
Finally, such that for every $k \in \N$, $u_k$ is the unique minimizer of
\begin{align*}
	w\longmapsto\int_\Omega \left[ F_k(Dw)+\e_k|Dw-Du_k| \right] \, dx,
\end{align*}
which is an integral functional with a Lipschitz continuous integrand. If we define 
\begin{equation*}
\sigma_k := F'_k(Du_k),
\end{equation*}
then we can derive a related important variational inequality. Indeed, fixed $\eta\in\K$ and $0<\e<1$, we have that defining
\begin{align*}
	v_k:=u_k+\e(\eta-u_k)\in\K
\end{align*}
we obtain that
\begin{align*}
	\frac{1}{\e}\int_\Omega \left[ F_k(Dv_k)+\e_k|Dv_k-Du_k|-F_k(Du_k) \right] \, dx\ge 0.
\end{align*}
On the other hand,
\begin{align*}
	&\frac{1}{\e}\int_\Omega \left[ F_k(Dv_k)+\e_k|Dv_k-Du_k|-F_k(Du_k) \right] \, dx\\
	= \, &\frac{1}{\e}\int_\Omega \left[ F_k(Du_k+\e D(\eta-u_k))+\e\e_k|D\eta-Du_k|-F_k(Du_k) \right] \, dx,
\end{align*}
so
\begin{align*}
	\frac{1}{\e}\int_\Omega \left[ F_k(Du_k+\e D(\eta-u_k))-F_k(Du_k)\right] \, dx \ge  -\e_k\int_\Omega|D\eta-Du_k| \, dx.
\end{align*}
Since
\begin{align*}
	\frac{1}{\e}\int_\Omega \left[ F_k(Du_k+\e D(\eta-u_k))-F_k(Du_k) \right] \, dx = \int_\Omega\int_0^1 \braket{F_k'(Du_k+s\e D(\eta-u_k)), D(\eta-u_k)} \, dsdx,
\end{align*}
then, by the Dominated Convergence Theorem and the uniform boundedness of $(Du_k)_k$ in $L^1(\Omega)$, we get
\begin{align}
	\label{dis-var}
	\int_\Omega \braket{F_k'(Du_k), D(\eta-u_k)} \, dx \ge -\e_k\int_\Omega|D\eta-Du_k| \, dx \ge -\e_k\left( C+\int_\Omega |D\eta| \, dx\right).
\end{align}


\subsection{Passage to the limit}
\label{quat_due}

\noindent
Now we want to understand what is the asymptotic behavior of $\sigma_k$ as $k\to \infty$. Since $F_k \nearrow F$ pointwise and $F_{k+1} \ge F_k$ for all $k \in \N$, it follows in particular from Dini's Lemma that the convergence is locally uniform in $\xi$, so we can now prove that $\sigma_k = F'_k(Du_k) \longrightarrow F'(Du)$ locally uniformly. To that end, we consider $\xi\in\R^n$ and $(F'_k(\xi_k))_k$ where $\xi_k \to \xi$ as $k\to \infty$. Because difference-quotients of convex functions are increasing in the increment, we have for all $\eta \in \R^n$ and $0 < |t| \le 1$ that
\begin{eqnarray*}
|\braket{F'_k(\xi_k) - F'(\xi), \eta}|
&\le& \left| \frac{F_k(\xi_k + t \, \eta) - F_k(\xi_k) - \braket{F'(\xi), t \, \eta}}{t} \right | \\
&\le& | F_k(\xi_k + \eta) - F_k(\xi_k) - \braket{F'(\xi), \eta} |.
\end{eqnarray*}
Consequently, we get
\begin{equation*}
\limsup_{k \, \to \, +\infty} | \braket{F'_k(\xi_k) - F'(\xi), \eta} | \le | F(\xi + \eta) - F(\xi) - \braket{F'(\xi), \eta} |,
\end{equation*}
for all $\eta \in \R^n$. Hence, for all $0 < s \le 1$, we obtain that
\begin{equation*}
\limsup_{k \, \to \, +\infty} \frac{ | \braket{F'_k(\xi_k) - F'(\xi), s \, \eta} | }{s} \le \frac{ | F(\xi + s \, \eta) - F(\xi) - \braket{F'(\xi), s \,\eta} | }{s},
\end{equation*}
which is equivalent to
\begin{equation*}
\limsup_{k \, \to \, +\infty} | \braket{F'_k(\xi_k) - F'(\xi), \eta}| \le \left| \frac{F(\xi + s \, \eta) - F(\xi)}{s} - \braket{F'(\xi), \eta} \right|.
\end{equation*}
If we let $s \to 0$ and recall that $F$ is differentiable in $\xi$, we conclude that the left-hand side must vanish. This proves the local uniform convergence of derivatives, so it follows that $\sigma_k = F'_k(Du_k) \longrightarrow F'(Du)$ locally uniformly and, in particular, in measure on $\Omega$.
\\
We now point out that the following identity holds, namely
\begin{equation}
\label{rel-pre-ext}
\braket{\sigma_k, Du_k} = F^*_k(\sigma_k) + F_k(Du_k),
\end{equation}
and it has been deduced by \eqref{fineq2} using the definition of $\sigma_k$, recalling that $F_k^{**} = F_k$ and choosing $\xi = Du_k$. Moreover, the same identity holds for the function $F$, namely
\begin{equation}
	\label{extr}
	\braket{\sigma^*, Du} = F^*(\sigma^*) + F(Du),
\end{equation}
with $\sigma^* := F'(Du)$.


\subsection{The validity of \eqref{dimdim}}
\label{quat_tre}

\noindent
Since $u_0 \in \K$, we can use \eqref{dis-var} choosing $\eta = u_0$, thus getting
\begin{equation*}
\int_\Omega \braket{\sigma_k, Du_0 - Du_k} \, dx \ge -\e_k\left( C+\int_\Omega|Du_0| \, dx\right)  \qquad \forall \, k \in \N,
\end{equation*}
that is equivalent to
\begin{align*}
	\int_\Omega \braket{\sigma_k, Du_0} \, dx + \e_k\left( C+\int_\Omega|Du_0| \, dx\right)\ge  \int_\Omega \braket{\sigma_k, Du_k} \, dx   \qquad \forall \, k \in \N.
\end{align*}
Now, integrating \eqref{rel-pre-ext} over $\Omega$ and using the previous inequality, we have that, chosen $t > 1$,
\begin{eqnarray*}
&& \int_\Omega F^*_k(\sigma_k) \, dx \\
&=& \int_\Omega \braket{\sigma_k, Du_k} \, dx - \int_\Omega F_k(Du_k) \, dx \\
&\le& \int_\Omega \braket{\sigma_k, Du_0} \, dx +\e_k\left( C+\int_\Omega|Du_0| \, dx\right) - \int_\Omega F_k(Du_k) \, dx \\
&=& \frac{1}{t} \int_\Omega \braket{\sigma_k, t \, Du_0} \, dx +\e_k\left( C+\int_\Omega|Du_0| \, dx\right) - \int_\Omega F_k(Du_k) \, dx \\
&\le& \frac{1}{t} \int_\Omega F^*_k(\sigma_k) \, dx + \, \frac{1}{t} \int_\Omega F_k( t \, Du_0) \, dx +\e_k\left( C+\int_\Omega|Du_0| \, dx\right) - \int_\Omega F_k(Du_k) \, dx,
\end{eqnarray*}
where we also used \eqref{fineq}, exploiting the convexity and lower semicontinuity of $F_k$. Reabsorbing the first term in the right-hand side by the left-hand side and using the fact that $F_k(\xi) \le F(\xi)$ for all $\xi \in \R^n$ we obtain that
\begin{eqnarray}
&& \int_\Omega F^*_k(\sigma_k) \, dx \nonumber \\
&\le& \frac{1}{t-1} \int_\Omega F_k( t \, Du_0) \, dx - \frac{t}{t-1} \int_\Omega F_k(Du_k) \, dx +\frac{t \,\e_k}{t-1}\left(C+\int_\Omega|Du_0|\,dx\right) \nonumber \\
&\le& \frac{1}{t-1} \int_\Omega F( t \, Du_0) \, dx - \frac{t}{t-1} \int_\Omega F_k(Du_k) \, dx +\frac{t\,\e_k}{t-1}\left(C+\int_\Omega|Du_0|\,dx\right) \nonumber \\
&\le& C_1 \int_\Omega F(t \, Du_0) \, dx - C_2 \int_\Omega F(Du_0) \, dx +C_3\left(C+\int_\Omega|Du_0|\,dx\right) \label{non-fin},
\end{eqnarray}
where in the last line we used the fact that
\begin{align*}
	\int_\Omega F_k(Du_k)\,dx\le \int_\Omega F(Du)\,dx\le \int_\Omega F(Du_0)\,dx,
\end{align*}
and that since $(\e_k)_k$ is converging then it is bounded. Recalling that $F^*_k \searrow F^*$ and, by the hypotheses, that $F(t \, Du_0) \in L^1(\Omega)$ with $t > 1$, we also obtain that
\begin{equation*}
\int_\Omega F^*(\sigma_k) \, dx \le C_1 \int_\Omega |F(t \, Du_0)| \, dx + C_2 \int_\Omega F(Du_0) \, dx +C_3\left(C+\int_\Omega|Du_0|\,dx\right) < +\infty.
\end{equation*}
Since we already observed that $\sigma_k \longrightarrow F'(Du)$ locally uniformly, by the previous estimate and Fatou's Lemma we have that
\begin{equation*}
\int_\Omega F^*(F'(Du)) \, dx \le \liminf_{k \, \to \, +\infty} \int_\Omega F^*(\sigma_k) \, dx < +\infty.
\end{equation*}
Thus
\begin{equation*}
F^*(F'(Du)) \in L^1(\Omega).
\end{equation*}
Whence, by \eqref{extr}, we also have
\begin{equation*}
\braket{F'(Du), Du} \in L^1(\Omega),
\end{equation*}
since $F(Du) \in L^1(\Omega)$ by the definition of minimizer.


\subsection{The validity of the variational inequality}
\label{quat_quat}

\noindent
Now we want to prove the validity of the variational inequality. Since $F^*$ is superlinear at infinity by Lemma \ref{superlinear} and since $F^*_k \searrow F^*$, then there exists $\overline{\theta} \, : \, [ 0, +\infty ) \longrightarrow [ 0, +\infty )$ increasing and superlinear at infinity such that 
\begin{equation}
\label{teta-polare}
\overline{\theta}(|\xi|) \le F^*(\xi) \le F^*_k(\xi) \qquad \forall \, \xi \in \R^n.
\end{equation}
Moreover, using \eqref{non-fin} we have that
\begin{equation*}
\sup_{k \, \in \, \N} \int_\Omega \overline{\theta}(|\sigma_k|) \, dx \le \int_\Omega F^*_k(\sigma_k) \, dx  < +\infty,
\end{equation*}
thus we can use De La Vall\`{e}-Poussin's Theorem in order to obtain the equi-integrability for $(\sigma_k)_k$. Since $(\sigma_k)_k$ converges in measure to $\sigma^*$, then we can apply Vitali's Convergence Theorem which proves that $(\sigma_k)_k$ converges to $\sigma^*$ in $L^1(\Omega)$. Now we observe that \eqref{dis-var} implies that
\begin{align}
	\label{dis-vark2}
	\int_\Omega\braket{\sigma_k,D\varphi} \, dx \ge -\e_k\int_\Omega|D\varphi| \, dx, \quad \forall\,\varphi\in C^\infty_0(\Omega), \varphi \ge 0.
\end{align}
Indeed this is implied by taking $\eta=u_k+\varphi$ in \eqref{dis-var}. It follows that if we pass to the limit as $k \longrightarrow +\infty$ in \eqref{dis-vark2} and, by the $L^1(\Omega)$ convergence of $(\sigma_k)_k$ to $\sigma^*$, it yields
\begin{equation*}
\int_\Omega \braket{\sigma^*, D\varphi} \, dx \ge 0 \qquad \forall \, \varphi \in \C^\infty_0(\Omega), \,\,\, \varphi \ge 0.
\end{equation*}
This implies that $\dive \, \sigma^* \le 0$ in the distributional sense, i.e. \eqref{div_neg}.


\subsection{The validity of \eqref{min_max_finale}}
\label{quat_cin}

\noindent
We should notice that up to subsequences it holds true that
\begin{equation}
\label{lim-inf}
\int_\Omega \braket{\sigma^*, Du} \, dx \le \liminf_{k \, \to \, +\infty} \int_\Omega \braket{\sigma_k, Du_k} \, dx.
\end{equation}
To prove this, we observe that by \eqref{Fkpositive}, \eqref{stime-acca} and \eqref{effe-kappa} we have
\begin{align*}
	F_k(\xi)\ge-\mu_k\qquad \forall\,\xi\in\R^n.
\end{align*}
Now, by \eqref{rel-pre-ext}, \eqref{teta-polare} and the fact that $\overline{\theta}(|\xi|) \ge 0$ for all $\xi \in \R^n$, we get that 
\begin{equation*}
\braket{\sigma_k, Du_k} = F^*_k(\sigma_k) + F_k(Du_k) \ge \overline{\theta}(|\sigma_k|) - \mu_k \ge - \mu_k.
\end{equation*}
Therefore we can apply Fatou's Lemma to the sequence of functions $(\braket{\sigma_k, Du_k} + \mu_k)_k$ that, thanks to the definition of $(\mu_k)_k$, converges a.e. to $\braket{\sigma^*, Du}$ up to subsequences, letting us deduce that
\begin{equation*}
\int_\Omega \braket{\sigma^*, Du} \, dx \le \liminf_{k \, \to \, +\infty} \int_\Omega \left(\braket{\sigma_k, Du_k} + \mu_k \right) \, dx \le \liminf_{k \, \to \, +\infty} \int_\Omega \braket{\sigma_k, Du_k} \, dx,
\end{equation*}
as we wanted to prove. Using \eqref{dis-var} with $\psi$ in place of $\eta$, since $\psi \in \K$, we get that
\begin{equation*}
\int_\Omega \braket{\sigma_k, Du_k} \, dx \le \int_\Omega \braket{\sigma_k, D\psi} \, dx +\e_k\left( C+\int_\Omega|D\psi|\,dx\right)
\end{equation*}
and the $L^1(\Omega)$ convergence of $(\sigma_k)_k$ to $\sigma^*$ implies that
\begin{align*}
&\liminf_{k \, \to \, +\infty} \int_\Omega \braket{\sigma_k, Du_k} \, dx \\
&\le \liminf_{k \, \to \, +\infty} \left[ \int_\Omega \braket{\sigma_k, D\psi} \, dx +\e_k\left( C+\int_\Omega|D\psi|\,dx\right)\right] \\
&= \int_\Omega \braket{\sigma^*, D\psi} \, dx.
\end{align*}
Combining this with \eqref{extr} and \eqref{lim-inf}, we obtain that
\begin{equation*}
\int_\Omega F(Du) \, dx + \int_\Omega F^*(\sigma^*) \, dx = \int_\Omega \braket{\sigma^*, Du} \, dx \le \int_\Omega \braket{\sigma^*, D\psi} \, dx = [\![\sigma^*, \psi]\!]_{u_0}(\overline{\Omega}),
\end{equation*}
where we use the measure defined in \eqref{misura-sigma}. Now we have to prove the reverse inequality, i.e.
\begin{equation}
\label{extr-min-max}
\int_\Omega F(Du) \, dx + \int_\Omega F^*(\sigma^*) \, dx \ge [\![\sigma^*, \psi]\!]_{u_0}(\overline{\Omega}).
\end{equation}
First we observe that, thanks to \eqref{parts},
\begin{eqnarray*}
[\![\sigma^*, \psi]\!]_{u_0}(\overline{\Omega})
&=& \int_\Omega (\psi - u_0) \, d(-\dive \, \sigma^*) + \int_\Omega \braket{\sigma^*, Du_0} \, dx \\
&=& \int_\Omega (\psi - u + u - u_0) \, d(-\dive \, \sigma^*) + \int_\Omega \braket{\sigma^*, Du_0} \, dx \\
&\le& \int_\Omega (u - u_0) \, d(-\dive \, \sigma^*) + \int_\Omega \braket{\sigma^*, Du_0} \, dx,
\end{eqnarray*}
where we used the fact that
\begin{equation*}
\int_\Omega (u - \psi) \, d(-\dive \, \sigma^*) \ge 0
\end{equation*}
because $-\dive \, \sigma^*$ is a non-negative Radon measure and $u \ge \psi$ a.e. on $\Omega$. By \eqref{parts} and by \eqref{extr}, we get that
\begin{eqnarray}
[\![\sigma^*, \, \psi]\!]_{u_0}(\overline{\Omega})
&\le&\int_\Omega (u - u_0) \, d(-\dive \, \sigma^*) + \int_\Omega \braket{\sigma^*, Du_0} \, dx \nonumber \\
&=&\int_\Omega \braket{\sigma^*, Du} \, dx \nonumber \\
&=&\int_\Omega F(Du) \, dx + \int_\Omega F^*(\sigma^*) \, dx. \label{misura-div2}
\end{eqnarray}
Combining \eqref{extr-min-max} with the previous inequality, we get that it is in fact an equality, i.e.
\begin{equation*}
[\![\sigma^*, \, \psi]\!]_{u_0}(\overline{\Omega}) = \int_\Omega F(Du) \, dx + \int_\Omega F^*(\sigma^*) \, dx.
\end{equation*}
Notice that this also implies that the measure $[\![\sigma^*, \, \psi]\!]_{u_0}(\overline{\Omega})$ is finite for $\sigma^*$ and $\psi$.

\section{Proof of Theorem \ref{dualproblem}}
\label{dimostr-2}

\noindent
The proof of Theorem \ref{dualproblem} leans on a similar method to the one used in Subsection \ref{quat_cin}, in particular exploiting the definition \eqref{misura-sigma} and the techniques of Convex Analysis. We write the details for the reader's convenience.
\\
	First of all, let us recall that since $u$ is a solution to \eqref{gprob} and $G$ is Lipschitz continuous, then it holds that
	\begin{equation}
		\label{varineq}
		\int_{\Omega}\braket{G'(Du), D\phi-Du} dx \geq 0 \qquad \forall \phi \in \K.
	\end{equation}
	Consider $\sigma \in S_{-}(\Omega)$ and $v \in \K$. Since $-\dive \sigma$ is a non-negative Radon measure and $v \ge \psi$ a.e. on $\Omega$, it holds that
	\begin{equation*}
		\int_\Omega (v - \psi) \, d(-\dive \sigma) \ge 0.
	\end{equation*}
	By definition \eqref{misura-sigma} and by the previous inequality we get that
	\begin{align*}
		[\![\sigma, \psi]\!]_{u_0}(\overline{\Omega})
		& =	\int_\Omega (\psi - u_0) \, d(-\dive \sigma) + \int_\Omega \braket{\sigma, Du_0} \, dx \\
		& = \int_\Omega (\psi - v + v - u_0) \, d(-\dive \sigma) + \int_\Omega \braket{\sigma, Du_0} \, dx \\
		& \le \int_\Omega (v - u_0) \, d(-\dive \sigma) + \int_\Omega \braket{\sigma, Du_0} \, dx.
	\end{align*}
	Since $v,u_0\in W_{u_0}^{1,1}(\Omega)$, then $v=u_0$ on $\partial \Omega$ and using \eqref{parts} on the first integral of the last expression, we get that
	\begin{align*}
		[\![\sigma, \psi]\!]_{u_0}(\overline{\Omega})
		& \le \int_\Omega \braket{\sigma, Dv - Du_0} \, dx + \int_\Omega \braket{\sigma, Du_0} \, dx \\
		& = \int_\Omega \braket{\sigma, Dv} \, dx \\
		& \le \int_\Omega G(Dv) \, dx + \int_\Omega G^*(\sigma) \, dx,
	\end{align*}
	by means of \eqref{fineq} and the fact that $G=G^{**}$. The previous inequality implies that
	\begin{equation*}
		\int_\Omega G(Dv) \, dx \ge [\![\sigma, \psi]\!]_{u_0}(\overline{\Omega}) - \int_\Omega G^*(\sigma) \, dx,
	\end{equation*}
	so passing to the minimum on $v \in \K$ and to the maximum on $\sigma \in S_{-}(\Omega)$ we get
	\begin{equation*}
		\min_{v \in \K} \int_\Omega G(Dv) \, dx \ge \max_{\sigma \in S_{-}(\Omega)} \left ([\![\sigma, \psi]\!]_{u_0}(\overline{\Omega}) - \int_\Omega G^*(\sigma) \, dx \right ).
	\end{equation*}
	We have now to prove the reverse inequality. We consider \eqref{fineq2} with $F=G$ and $\xi=Du$ and, since $u \in \K$ is the unique solution of \eqref{gprob}, then
	\begin{align*}
		\int_\Omega G(Du) \, dx
		& = \int_\Omega \braket{G'(Du), Du} \, dx - \int_\Omega G^*(G'(Du)) \, dx \\
		& = \int_\Omega \braket{G'(Du), Du - Du_0} \, dx + \int_\Omega \braket{G'(Du), Du_0} - \int_\Omega G^*(G'(Du)) \, dx.
	\end{align*}
	Since $G$ satisfies the variational inequality \eqref{varineq} then it holds true also that
	\begin{equation*}
		\int_\Omega \braket{G'(Du), D\phi} \, dx \ge 0 \qquad \forall \phi \in C_0^\infty(\Omega), \phi \ge 0,
	\end{equation*}
	and so if we set $\sigma = G'(Du)$ then we have that $\sigma \in S_{-}(\Omega)$, thus we have
	\begin{eqnarray*}
		\int_\Omega G(Du) \, dx
		&=& \int_\Omega \langle \sigma, Du - Du_0 \rangle \, dx + \int_\Omega \langle \sigma, Du_0 \rangle \, dx - \int_\Omega G^*(\sigma) \, dx \\
		&=& \int_\Omega (u - u_0) \, d (-\dive \sigma) + \int_\Omega \langle \sigma, Du_0 \rangle \, dx - \int_\Omega G^*(\sigma) \, dx \\
		&=& \int_\Omega (\psi - \psi + u - u_0) \, d (-\dive \sigma) + \int_\Omega \langle \sigma, Du_0 \rangle \, dx - \int_\Omega G^*(\sigma) \, dx \\
		&=& \int_\Omega (\psi - u_0) \, d (-\dive \sigma) - \int_\Omega (\psi - u) \, d (-\dive \sigma) + \int_\Omega \langle \sigma, Du_0 \rangle \, dx - \int_\Omega G^*(\sigma) \, dx \\
		&=& [\![\sigma, \psi]\!]_{u_0}(\overline{\Omega}) - \int_\Omega \langle \sigma, D\psi - Du \rangle \, dx - \int_\Omega G^*(\sigma) \, dx \\
		&\le& [\![\sigma, \psi]\!]_{u_0}(\overline{\Omega}) - \int_\Omega G^*(\sigma) \, dx,
	\end{eqnarray*}
	where we used once again \eqref{parts} and the fact that it holds \eqref{varineq}, where we chose $\phi = \psi \in \K$.
	Finally, we have that
	\begin{align*}
		\min_{v \in \K} \int_\Omega G(Dv) \, dx 
		&= \int_\Omega G(Du) \, dx \\
		& \le [\![\sigma, \psi]\!]_{u_0}(\overline{\Omega}) - \int_\Omega G^*(\sigma) \, dx \\
		& \le \max_{\sigma \in S_{-}(\Omega)} \left([\![\sigma, \psi]\!]_{u_0}(\overline{\Omega}) - \int_\Omega G^*(\sigma) \, dx \right),
	\end{align*}
Since this also implies that the measure $[\![\sigma, \, \psi]\!]_{u_0}(\overline{\Omega})$ is finite for $\sigma$ and $\psi$, this concludes the proof.

\color{black}


\section*{Acknowledgments}
\noindent
S.R. acknowledges financial support from the Austrian Science Fund (FWF) projects F65 and Y1292.


\end{document}